\documentclass[12pt]{amsart}
\usepackage{amsmath,amsthm,amsfonts,amscd,amssymb,eucal,latexsym,mathrsfs,enumerate,framed}

\usepackage{hyperref}

\newcommand{\Bb}{\mathbb B}
\newcommand{\Hs}{\mathcal H}

\newcommand{\ad}{\mathrm{ad}}
\newcommand{\eps}{\varepsilon}

\newcommand{\A}{\mathcal A}
\newcommand{\B}{\mathcal B}
\newcommand{\M}{\mathcal M}
\newcommand{\N}{\mathcal N}

\theoremstyle{plain}

\newtheorem*{theorem*}{Theorem}

\newtheorem*{proposition*}{Proposition}

\theoremstyle{definition}
\newtheorem*{definition*}{Definition}

\theoremstyle{remark}

\title{A remark on the similarity and perturbation problems}
\author[Cameron]{Jan Cameron}
\address{\hskip-\parindent
Jan Cameron, Department of Mathematics, Vassar College, Poughkeepsie, NY 12604,
U.S.A.}
\email{jacameron@vassar.edu}
\thanks{JC's research is partially supported by an AMS-Simons research travel grant.}
\author[Christensen]{Erik Christensen}
\address{\hskip-\parindent
Erik Christensen, Institute for Mathematiske Fag, University of Copenhagen, Copenhagen, Denmark.}
\email{echris@math.ku.dk}
\author[Sinclair]{Allan M.~Sinclair}
\address{\hskip-\parindent
Allan M.~Sinclair, School of Mathematics, University of Edinburgh, JCMB, King's Buildings, Mayfield Road, Edinburgh, EH9 3JZ, Scotland.}
\email{a.sinclair@ed.ac.uk}
\author[Smith]{Roger R.~Smith}
\address{\hskip-\parindent
Roger R.~Smith, Department of Mathematics, Texas A{\&}M University,
College Station, TX 77843, U.S.A.}
\email{rsmith@math.tamu.edu}
\thanks{RS is partially supported by NSF grant DMS-1101403}

\author[White]{Stuart White}
\address{\hskip-\parindent
Stuart White, School of Mathematics and Statistics, University of Glasgow, 
University Gardens, Glasgow Q12 8QW, Scotland.}
\email{stuart.white@glasgow.ac.uk}
\thanks{SW is partially supported by EPSRC grant EP/I019227/1}

\author[Wiggins]{Alan D.~Wiggins}
\address{\hskip-\parindent
Alan D.~Wiggins, Department of Mathematics and Statistics. University of Michigan-Dearborn, 4901 Evergreen Road, Dearborn, MI 48126, USA.}
\email{adwiggin@umd.umich.edu}

\date{\today}

\begin{document}

\begin{abstract}
In this note we show that Kadison's similarity problem for $C^*$-algebras is equivalent to a problem in perturbation theory: must close $C^*$-algebras have close commutants?\\

\noindent{\sc R\'esum\'e.} 
Dans cet note, nous montrons que le probl\`e{}me de similarit\'e{} de Kadison est \'e{}quivalent \`a{} la question suivante en th\'e{}orie de la perturbation: les commutants de deux $C^*$-alg\`e{}bres proches sont-ils n\'e{}cessairement proches?
\end{abstract}

\maketitle

Let $\A$ be a $C^*$-algebra. In 1955 Kadison asked whether every bounded homomorphism $\sigma:\A\rightarrow\Bb(\Hs)$ is similar to a
$^*$-homomorphism \cite{K:AJM}, i.e. does there exist an invertible operator
$S\in\Bb(\Hs)$ such that $S^{-1}\sigma(\cdot)S$ is a $^*$-homomorphism? This
problem, which remains open, has positive answers in the following cases: for
amenable algebras \cite{B:PAMS}, for properly infinite von Neumann algebras and
traceless $C^*$-algebras \cite{H:Ann}, in the presence of a cyclic vector
\cite{H:Ann},  and for II$_1$ factors with property $\Gamma$ \cite{C:JOT}. It is also
equivalent to a number of important problems and properties:
\begin{itemize}
\item is every bounded homomorphism $\sigma : \A\rightarrow\Bb(\Hs)$ completely
bounded? (\cite{H:Ann})
\item the \emph{derivation problem} \cite{Ki:JOT} - given a faithful non-degenerate representation $\iota:\A\rightarrow\Bb(\Hs)$, is every derivation of $\iota(\A)$ into $\Bb(\Hs)$ spatial?
\item the \emph{distance property} \cite{C:JFA,C:Scand,Ki:JOT} - does there exist a constant $K>0$ such that given any faithful non-degenerate representation $\iota:\A\rightarrow\Bb(\Hs)$ 
\begin{equation}\label{0}
\frac{1}{2}\|\ad(T)|_{\iota(\A)}\|_{cb}=d(T,\iota(\A)')\leq K\|\ad(T)|_{\iota(\A)}\|
\end{equation}
holds for all $T\in\Bb(\Hs)$?\footnote{The equality in (\ref{0}) is Arveson's distance formula from \cite{A:JFA} and is always valid.} When this holds, $\A$ is said to have property $D_K$. Here, and throughout the paper, $\ad(T)$ denotes the spatial derivation $X\mapsto [T,X]=TX-XT$ on $\Bb(\Hs)$.
\item do all $C^*$-algebras have finite length?\footnote{Pisier's
characterisation of the similarity problem using his notion of {\em{length}}
provides an intrinsic formulation of the similarity property in terms of
internal behaviour.} (\cite{Pi:StP})
\item are all von Neumann algebras hyperreflexive? (see \cite[Theorem 10.3]{P:Book}).
\item Diximer's invariant operator range problem from 1950 (see \cite[Problem 10.4]{P:Book}).
\end{itemize}
In this note we add to the list of equivalent formulations by showing that the similarity problem is  equivalent to an open question in the uniform perturbation theory of operator algebras.

The Kadison-Kastler metric on all $C^*$-subalgebras of $\Bb(\Hs)$ is given by
restricting the Hausdorff metric on subsets of $\Bb(\Hs)$ to the unit balls.
Specifically, for pairs of $C^*$-algebras $\A,\B\subseteq\Bb(\Hs)$,
$$
d(\A,\B)=\max\left\{\sup_{\substack{x\in \A\\ \|x\|\leq
1}}\,\inf_{\substack{y\in
\B\\ \|y\|\leq 1}}\|x-y\|,\,\sup_{\substack{y\in \B\\ \|y\|\leq
1}}\,\inf_{\substack{x\in \A\\ \|x\|\leq 1}}\|y-x\|\right\}.
$$
This notion was introduced in \cite{KK:AJM} which conjectures that sufficiently
close von Neumann algebras must arise from a small unitary perturbation. Precisely, for all $\eps>0$, there should exist $\delta>0$ such that, given von Neumann
algebras $\M,\N\subseteq\Bb(\Hs)$ with $d(\M,\N)<\delta$, there exists a unitary
$U$ on $\Hs$ with $\|U-I_{\Hs}\|<\eps$ satisfying $U\M U^*=\N$.\footnote{Due to
the results of \cite{CC:BLMS,J:CMB}, the analogous conjecture for $C^*$-algebras
is that sufficiently close separable $C^*$-algebras acting on a separable
Hilbert space should be spatially isomorphic; one cannot demand control on the
distance from an implementing unitary to the identity in this context.}  This
conjecture, if true, immediately implies that the operation of taking commutants
of $C^*$-algebras is continuous with respect to the Kadison-Kastler
metric.\footnote{Via a Kaplansky density argument, it suffices to consider the
continuity of taking commutants on the set of von Neumann subalgebras of
$\Bb(\Hs)$.} Establishing that   certain close
operator algebras have close commutants,
and a related property of near inclusions (see Remark 2 below) play a key
role in our recent work on the Kadison-Kastler problem
 for crossed product factors \cite{CCSSWW:Preprint}, in EC's solution to the
near inclusion problem for injective von Neumann algebras \cite{C:Acta} and in
work on transferring $K$-theoretic invariants between close algebras
\cite{Kh:JOT,CSSW:GAFA,PTWW:InPrep} (see Remark 3 below).  Here we show that
the continuity of the operation of taking commutants is equivalent to the
similarity problem. In order to phrase this result for individual $C^*$-algebras
we make the following definition.

\begin{definition*}
Let $\A$ be a $C^*$-algebra. Say that \emph{commutants are continuous at $\A$}
if, for $\eps>0$, there exists $\delta>0$ such that given a faithful
non-degenerate $^*$-representation $\iota:\A\rightarrow\Bb(\Hs)$ and a
$C^*$-algebra $\B\subseteq\Bb(\Hs)$ with $d(\iota(\A),\B)<\delta$, we have
$d(\iota(\A)',\B')<\eps$.
\end{definition*}

\begin{theorem*}
The following statements are equivalent for a $C^*$-algebra $\A$:
\begin{enumerate}[(I)]
\item\label{T.1} $\A$ satisfies Kadison's similarity property;
\item\label{T.2} commutants are continuous at $\A$.
\end{enumerate}
Further, the similarity problem has a positive answer for all $C^*$-algebras if
and only if the operation of taking commutants on $\Bb(\Hs)$ is continuous
in the case when $\Hs$ is an infinite dimensional separable Hilbert
space.
\end{theorem*}

\section*{Proof of the Theorem}
If the similarity problem has a positive answer for all $C^*$-algebras, then there exists a constant $K>0$ such that all $C^*$-algebras have property $D_K$ (by \cite{Ki:JOT} and \cite{C:Scand}). It is then easy to see that
\begin{equation}\label{NewEq}
d(\A',\B')\leq 4K d(\A,\B)
\end{equation}
for all $C^*$-algebras $\A,\B\subset\Bb(\Hs)$ (this dates back to \cite{C:Acta}, and can be found explicitly by combining Proposition 2.3 (iii) and Proposition 2.5 of \cite{CSSW:GAFA}).  The local version (\ref{T.1})$\implies$(\ref{T.2}) of this implication is much more involved as we do not have the information that $\B$ also has property $D_K$ required in the previous argument. Nevertheless, this implication was established in \cite[Theorem 4.2]{CSSW:GAFA} as a step towards showing that if $\A$ has the similarity property, then any $C^*$-algebra sufficiently close to $\A$ also has the similarity property.

We now turn to the reverse implication (\ref{T.2})$\implies$(\ref{T.1}). For $X,R\in\Bb(\Hs)$,  use the identity
$R^mX-XR^m=R^{m-1}(RX-XR)+(R^{m-1}X-XR^{m-1})R$ to establish by induction the
estimates
\begin{equation}\label{1}
\|[R^m,X]\|\leq m\|R\|^{m-1}\|[R,X]\|,\quad \quad m\in\mathbb N.
\end{equation}
 
Suppose that commutants are continuous at $\A$, but that there
exists no constant $K>0$ such that $\A$ has property $D_K$.   Then we may find
a sequence $(\iota_n)_n$ of faithful, non-degenerate representations
$\iota_n:\A\rightarrow\Bb(\Hs_n)$ and operators $T_n\in\Bb(\Hs_n)$ with
$$
d(T_n,\iota_n(\A)')>n\|\ad(T_n)|_{\iota_n(\A)}\|, \quad \quad n\geq 1.
$$
For each $n$, choose $\widetilde{S_n}\in\iota_n(\A)'$ satisfying
$\|\widetilde{S_n}-T_n\|=d(T_n,\iota_n(\A)')$. Write $S_n=T_n-\widetilde{S_n}$ so
that 
\begin{equation}\label{4}
\|S_n\|>n\|\ad(S_n)|_{\iota_n(\A)}\|,\quad \quad n\geq 1.
\end{equation}
Scaling in this inequality allows us to assume that $\|S_n\|=1$ while retaining
the validity of  (\ref{4}). Since $d(S_n,\iota_n(\A)')=1$, either
$d(\Re(S_n),\iota_n(\A)')\geq 1/2$ or $d(\Im(S_n),\iota_n(\A)')\geq 1/2$ for
each $n$, while the equation
$\|\ad(S_n)|_{\iota_n(\A)}\|=\|\ad(S_n^*)|_{\iota_n(\A)}\|$ gives the limiting
behaviour
$\|\ad(\Re(S_n))|_{\iota_n(\A)}\|,\|\ad(\Im(S_n))|_{\iota_n(\A)}\|\rightarrow 0$
as $n\rightarrow\infty$.  Thus we obtain a sequence of self-adjoint contractions
 $R_n\in\Bb(\Hs_n)$ such that $d(R_n,\iota_n(\A)')\geq
1/2$ for $n\geq 1$ and $\lim_{n\to\infty}\|\ad(R_n)|_{\iota_n(\A)}\|=
0$. By approximating $f(t)=\sqrt{1-t^2}$ uniformly by polynomials on
$[-1,1]$ and using the estimate (\ref{1}), we see that  the unitary
operators $$
U_n^{\pm}=R_n\pm i\sqrt{1-R_n^2}
$$
satisfy $\lim_{n\to\infty}\|\ad(U_n^{\pm})|_{\iota_n(\A)}\|= 0$. Since
$R_n=(U_n^++U_n^-)/2$ for each $n$, either $d(U_n^+,\iota_n(\A)')\geq 1/2$ or
$d(U_n^-,\iota_n(\A)')\geq 1/2$.  Appropriate choices of signs give unitaries
$U_n\in\Bb(\Hs_n)$ with
\begin{equation}\label{5}
d(U_n,\iota_n(\A)')\geq 1/2,\quad\quad n\geq 1,
\end{equation}
while
$\lim_{n\to\infty}\|\ad(U_n)|_{\iota_n(\A)}\|= 0$. 

Now consider the faithful non-degenerate representations
$\iota_n\oplus\iota_n:\A\rightarrow {\mathbb{M}}_2(\Bb(\Hs_n))$ so that 
$$
(\iota_n\oplus\iota_n)(\A)=\left\{\begin{pmatrix}\iota_n(a)&0\\0&\iota_n(a)\end{pmatrix}:a\in \A\right\}.
$$
Define $C^*$-algebras $\B_n$ by
$$
\B_n=\left\{\begin{pmatrix}\iota_n(a)&0\\0&U_n\iota_n(a)U_n^*\end{pmatrix}:a\in
\A\right\}\subseteq{\mathbb{M}}_2(\Bb(\Hs_n)),\quad\quad n\geq 1,
$$
so that $d((\iota_n\oplus\iota_n)(\A),\B_n)\leq
\|\ad(U_n)|_{\iota_n(\A)}\|\rightarrow 0$.  Clearly
$W_n=\begin{pmatrix}0&I_{\Hs_n}\\0&0\end{pmatrix}\in
(\iota_n\oplus\iota_n)(\A)'$ so,
by hypothesis, $d(W_n,\B_n')\rightarrow 0$. Choose $\widetilde{W_n}\in\B_n'$
with
$\|W_n-\widetilde{W_n}\|=d(W_n,\B_n')$ and let  $X_n$ be the $(1,2)$-entry of
$\widetilde{W_n}$. Then $\|X_n-I_{\Hs_n}\|\rightarrow 0$ and
$$
\iota_n(a)X_n=X_nU_n\iota_n(a)U_n^*,\quad a\in \A.
$$
Thus $X_nU_n\in \iota_n(\A)'$ so that $d(U_n,\iota_n(\A)')\rightarrow 0$,
contradicting (\ref{5}). This  completes the proof of
(\ref{T.2})$\implies$(\ref{T.1}).

Finally, note that Kadison's similarity problem is equivalent to the existence of a
constant $K>0$ such that 
$$
\frac{1}{2}\|\ad(T)|_\A\|_{cb}\leq K\|\ad(T)|_\A\|,\quad T\in\Bb(\Hs),
$$
whenever $\Hs$ is a \emph{separable} Hilbert space and $\A\subseteq\mathbb B(\Hs)$ is a
$C^*$-algebra. This is established by a standard argument (cf the last two paragraphs of the proof of
\cite[Lemma 2.1.5]{CCSSWW:Preprint}) which shows that this last condition
implies the same condition without the separability assumption. Consequently, the proof of (\ref{T.2})$\implies$(\ref{T.1}) shows that if the operation of taking commutants is continuous on $\Bb(\Hs)$ for $\Hs$ separable and infinite dimensional, then Kadison's similarity property holds.\qed
\section*{Remarks}
\subsection*{1}
If the operation of taking commutants is globally continuous, then there is a universal constant $K>0$ such that all $C^*$-algebras have property $D_K$ and hence (\ref{NewEq}) holds.  In particular taking commutants is Lipschitz.  In the local situation, the proof of (\ref{T.1})$\implies$(\ref{T.2}) given in \cite[Theorem
4.2]{CSSW:GAFA} shows that if $\A\subseteq\Bb(\Hs)$ is a $C^*$-algebra which has
finite similarity length $\ell$ and length constant $K$, then there exists a
constant $C'_{\ell,K}$ such that 
$$
d(\A',\B')\leq C'_{\ell,K} d(\A,\B)
$$
for all $C^*$-algebras $\B$ on $\Hs$ with $d(\A,\B)$ sufficiently small. Since the metric $d$ has diameter $1$, it follows that there is a constant $C_{\ell,K}$ such that
$$
d(\A',\B')\leq C_{\ell,K} d(\A,\B)
$$
for all $C^*$-algebras $\B$ on $\Hs$. This says that if $\A$ has the
similarity property, then taking commutants is locally a Lipschitz operation
near $\A$.  

\subsection*{2}  The similarity property is also characterised by the property
that near inclusions of algebras give rise to near inclusions of their
commutants. For $C^*$-algebras $\A,\B\subseteq\Bb(\Hs)$ and $\gamma\geq 0$,
write $\A\subseteq_\gamma\B$ if, given $x\in\A$,  there exists $y\in\B$ with
$\|x-y\|\leq\gamma\|x\|$.  It is natural to ask whether, given $\eps>0$, there
exists $\delta>0$ such that the near inclusion $\A\subseteq_\delta\B$ implies
 that $\B'\subseteq_\eps\A'$.\footnote{One can
also formulate this question locally at $\A$ by asking that this hold under all
faithful non-degenerate representations of $\A$.} Indeed the connection between
the similarity problem and questions of close commutants originates in
\cite{C:Acta} where a positive answer is given when $\A$ has property
$D_K$.  This question is also equivalent to the similarity property, and the
local version at $\A$ is equivalent to the similarity property for $\A$.   The
global statement follows immediately, as the ability to take commutants of all
near inclusions in this way implies that taking commutants is continuous. For
the local statement, note that the proof of the theorem only uses the continuity
of commutants for pairs $(\A,U\A U^*)$, and this follows from being able to take
commutants of near inclusions of $\A$.

\subsection*{3} It is natural to consider a completely bounded version of the
Kadison-Kastler metric: $d_{cb}(\A,\B)=\sup_nd(\A\otimes {\mathbb{M}}_n,B\otimes
{\mathbb{M}}_n)$.\footnote{Similarly define $\A\subset_{cb,\gamma}\B$ if and
only if $\A\otimes
{\mathbb{M}}_n\subset_\gamma\B\otimes {\mathbb{M}}_n$ for all $n\geq 1$.} When 
$d_{cb}(\A,\B)$ is
sufficiently small, one can show that $\A$ and $\B$ have the same $K$-theoretic
invariants used in the classification programme, \cite{Kh:JOT,CSSW:GAFA}.
Moreover in this case $\A$ and $\B$ also have isomorphic Cuntz semigroups \cite{PTWW:InPrep}. Thus
it is of interest to learn when $d(\cdot,\cdot)$ and $d_{cb}(\cdot,\cdot)$ are
equivalent metrics: this too is a reformulation of the similarity
problem.\footnote{A local version of this statement also holds: $\A$ has the
similarity property if and only if $d_{cb}(\A,\cdot)$ is equivalent to
$d(\A,\cdot)$.} That a positive solution to 
the similarity problem implies the equivalence of these metrics is noted at the
end of \cite[Section 4]{CSSW:GAFA}. The converse follows as Arveson's distance
formula shows that $\A\subset_{cb,\gamma}\B\implies \B'\subset_{cb,\gamma}\A'$
(see \cite[Proposition 2.2.3]{CCSSWW:Preprint}) and so $$d(\A',\B')\leq
d_{cb}(\A',\B')\leq 2d_{cb}(\A,\B).$$

\subsection*{Acknowledgements} When SW presented the work of \cite{CCSSWW:Preprint} in Vanderbilt University in January 2012, Dietmar Bisch asked whether it was necessary to use the similarity property in order to take commutants of close $C^*$-algebras. This directly sparked the present paper and we are grateful to Dietmar for bringing  this question to our attention.

\providecommand{\bysame}{\leavevmode\hbox to3em{\hrulefill}\thinspace}
\providecommand{\MR}{\relax\ifhmode\unskip\space\fi MR }
\providecommand{\MRhref}[2]{%
  \href{http://www.ams.org/mathscinet-getitem?mr=#1}{#2}
}
\providecommand{\href}[2]{#2}

\end{document}